\documentclass[letterpaper, 10 pt, conference]{ieeeconf}  % Comment this line out
\IEEEoverridecommandlockouts                              % This command is only
\usepackage{graphics} % for pdf, bitmapped graphics files
\usepackage{epsfig} % for postscript graphics files% \usepackage{epstopdf}
\usepackage{amsmath} % assumes amsmath package installed
\usepackage{amssymb}  % assumes amsmath package installed
\usepackage{caption}

\usepackage{verbatim}
\usepackage{graphicx} % support the \includegraphics command and options
\usepackage{epsfig} % for postscript graphics files% \usepackage{epstopdf}
\usepackage{epstopdf}
\usepackage{amsmath} % assumes amsmath package installed
\usepackage{amssymb}  % assumes amsmath package installed
\usepackage{color}
\usepackage{tikz}
\usetikzlibrary{matrix}

\usepackage{lipsum}
\makeatletter
\def\footnoterule{\relax%
	\kern-5pt
	\hbox to \columnwidth{\hfill\vrule width .9\columnwidth height 0.4pt\hfill}
	\kern4.6pt}
\makeatother

\makeatletter
\newcommand{\figcaption}{\def\@captype{figure}\caption}
\newcommand{\tabcaption}{\def\@captype{table}\caption}

\usepackage[hidelinks]{hyperref}

\makeatother
\title{\LARGE \bf The Impact of Execution Delay on Kelly-Based Stock\\
	Trading: High-Frequency  Versus Buy and Hold
}

\author{\large Chung-Han Hsieh,$^{*}$  B. Ross Barmish,$^{**}$ and John A. Gubner$^{***}$% <-this % stops a space
	\thanks{\hskip -10pt ${}^*$Chung-Han Hsieh is a graduate student working towards a Ph.D. degree in the Department of Electrical and Computer Engineering, University of Wisconsin, Madison, WI 53706. E-mail: \href{mailto: hsieh23@wisc.edu}{hsieh23@wisc.edu}} 
	\thanks{\hskip -10pt ${}^{**}$B. Ross Barmish is a \, Research \, Professor \, in \, the \,
		Department \, of \, Elec\-tric\-al and Computer  Engineering, Boston University, Boston, MA\ \ 02215.
		E-mail: \href{mailto: barmish@bu.edu}{barmish@bu.edu} 	}
	\thanks{\hskip -10pt ${}^{***}$John A. Gubner is a Professor in  the Department of Electrical and Computer Engineering, University of Wisconsin, Madison, WI 53706. \mbox{E-mail}: \href{mailto: john.gubner@wisc.edu.}{john.gubner@wisc.edu} 	  \vspace{3mm} }
}

\begin{document}
	
	\maketitle
	\thispagestyle{empty}
	\pagestyle{empty}
	
	\parindent = 0pt
	%%%%%%%%%%%%%%%%%%%%%%%%%%%%%%%%%%%%%%%%%%%%%%%%%%%%%%%%%%%%%%%%%%%%%%%%%%%%%%%%
	
	\begin{abstract}
		Stock trading based on Kelly's celebrated Expected Logarithmic Growth (ELG) criterion, a well-known prescription for optimal resource allocation, has received considerable attention in the literature. Using ELG as the performance metric, we compare the impact of trade execution delay on the relative performance of high-frequency trading versus buy and hold. While it is intuitively obvious and straightforward to prove that in the presence of sufficiently high transaction costs, buy and hold is the better strategy, is it possible that with no transaction costs, buy and hold can still be the better strategy? When there is no delay in trade execution, we prove a theorem saying that the answer is ``no.'' However, when there is delay in trade execution, we present simulation results using a binary lattice stock model to show that the answer can be ``yes.'' This is seen to be true whether self-financing is imposed or not.
	\end{abstract}

	%%%%%%%%%%%%%%%%%%%%%%%%%%%%%%%%%%%%%%%%%%%%%%%%%%%%%%%%%%%%%%%%%%%%%%%%%%%%%%%%%%
	
	\vspace{8mm}
	\section{Introduction}
	\label{SEC: INTRODUCTION}
	\vspace{0mm}
	Stock trading based on Kelly's celebrated  \textit{Expected Logarithmic Growth} (ELG) criterion, a well-known prescription for optimal resource allocation, has received considerable attention in the literature.  
	The formulation, first introduced in a betting scenario in the  seminal paper~\cite{Kelly_1956}, has been extended to address stock trading and portfolio rebalancing problems; e.g., see~\mbox{\cite{Hakansson_1971}--\cite{Hsieh_Barmish_Gubner_2016_CDC}}. 
	The reader is also referred to~\cite{MacLean_Thorp_Ziemba_2011} for a rather comprehensive exposition covering many aspects of  the  theory. 
	
	\vspace{3mm}
	This paper is most closely related to more recent work such as~\mbox{\cite{Kuhn_Luenberger_2010}--\cite{Hsieh_Gubner_Barmish_2018_CDC}} which provide results on the effect of rebalancing frequency on optimal trading performance.
	Specifically, in~\cite{Hsieh_Barmish_Gubner_2018_ACC} and~\cite{Hsieh_Gubner_Barmish_2018_CDC}, it is shown that in a so-called idealized market with a stock satisfying  a certain ``sufficient attractiveness" condition,  the buy and holder can match the performance  of the high-frequency  trader. Additionally, in~\cite{Hsieh_Barmish_Gubner_2018_ACC}, it was shown that when transaction costs are added into the mix, consistent with intuition, the buy and holder can strictly outperform the high-frequency trader. In~\cite{Hsieh_Gubner_Barmish_2018_CDC}, the question is raised whether this strict out-performance can happen when there are no transaction costs. In this paper, we prove that it cannot. In other words, for this case, high-frequency trading is ``unbeatable" in terms of expected logarithmic growth.

	\vspace{3mm}
	This result brings us to the next question which we address in this paper: If there is a delay in the trading system, sometimes called \textit{latency}, is high-frequency trading still unbeatable? 
	%This result brings us to the second goal of this paper. 
	To this end, here we consider the practical issue of delay in trade execution, which has not been considered to date in the existing ELG literature.
	In this regard, we introduce a one-step delay in trade execution into our formulation and consider both the self-financed and leveraged cases. In this context, our goal is to raise the possibility that when such a delay is present in real-world financial markets,  the  buy-and-hold strategy can achieve strictly higher ELG performance than  high-frequency~trading.    
	We use the word ``possibility" here because our formal results are obtained using a mathematical model with returns $X(k)$ which are larger than those seen with real ``high-frequency" trading data. 
	In view of this issue, in the final section, based on a simulation using high-frequency historical tick data leading to suggestive, but numerically inconclusive~results, we suggest future directions for research.
	%To this end, we use a binary lattice stock model as example to demonstrate this possibility. 
	%	With this as backdrop, it is natural to ask whether the buy and holder can achieve a higher ELG than a high-frequency trader in the real world.
	
	% 
	%
	%
	%
	
	\vspace{3mm}
	The plan for the remainder of this paper is as follows: In Section~\ref{SEC: Problem Formulation}, for the sake of self-containment, we summarize  our frequency-based formulation introduced in~\cite{Hsieh_Barmish_Gubner_2018_ACC} and~\cite{Hsieh_Gubner_Barmish_2018_CDC}. 
	In Section~\ref{SEC: The Blessing Of High-Frequency Trading},  we consider the no-delay case, and we provide a new result which we call the \textit{High-Frequency Maximality Theorem}. This theorem indicates that in the absence of transaction costs, high-frequency trading is \textit{unbeatable} in the sense of expected logarithmic growth.
	Then in Section~\ref{SEC: Extension to Include Delay In Trade Execution}, we extend the formulation to include execution delay and self-financing considerations. 
	In Section~\ref{SEC: Example Include Execution Delay}, working with this new formulation, we use a binary lattice stock price model to demonstrate that  the buy-and-hold strategy can outperform high-frequency trading. 
	In Section~\ref{SEC: Conclusion and Future Work}, some concluding remarks are provided and, per discussion above, an approach to study trading performance using high-frequency trading historical data is given, and some possible directions for future research are~indicated. 

	\vspace{8mm}
	\section[Frequency-Based Problem Formulation]{Frequency-Based Problem Formulation }
	\label{SEC: Problem Formulation}
	\vspace{0mm}
	A single stock is considered for trading over a finite time window at prices
	$S(k)>0$ for $k=0,\ldots,n$.
	The \textit{time} between stage~$k$ and~\mbox{$k+1$}, call it $\Delta t$, is viewed as small in the spirit of high-frequency trading; i.e., a fraction of a second. Additionally, we assume that stock-trading occurs within an idealized market. 
	That is, we assume  zero transaction costs, zero interest rates and perfect liquidity conditions. There is no gap between the bid and ask prices, and  the trader can buy or sell any number of shares including fractions at the traded price $S(k)$.  
	For more details on the idealized market assumption, the reader is referred to reference~\cite{Barmish_Primbs_TAC_2016}.

	\vspace{3mm}
	In this setting, we compare the performance of two traders. The first is a high-frequency trader who submits an order at each stage,  and the other is a buy and holder who only submits one order at $k=0$. Let the trader's account value and investment in dollars at time $k$ be denoted by~$V(k)$ and~$I(k)$, respectively. 
	We require that all trades be long-only, that is,~$I(k)\ge 0$, and self-financed; i.e., \mbox{$I(k) \le V(k)$}. Together, these imply that we also  require $V(k)\ge 0$. Now in the Kelly framework, discussed later in this section, the trader's investment level is given by a linear feedback~\mbox{$I(k)=K V(k)$}; e.g., see~\cite{Hsieh_Barmish_Gubner_2016_CDC} and \cite{Hsieh_Barmish_Gubner_2018_ACC}--\cite{Rujeerapaiboon_Barmish_Kuhn_2018}.
	%Consistent with the literature,  we adopt the same the linear feedback.
	Hence, the long-only condition leads to the requirement $K\ge 0$.
	% $I(k) \geq 0$. e self-financing condition
	%, that the investment level cannot exceed the account value, means~$I(k) \leq V(k)$ for all $k$. This
	The self-financing requirement corresponds to $K \leq 1$ when there is no execution delay. On the other hand, when execution delay is present, as seen in Section~\ref{SEC: Extension to Include Delay In Trade Execution}, this leads to constraint~\mbox{$K \leq {1}/(1+X_{\max})$}.

	\vspace{3mm}
	In the sequel, we primarily work with the \textit{returns}
	\[
	X(k) \doteq \frac{S(k+1) - S(k)}{S(k)}
	\]
	for $k=0,\ldots,n-1$.
	We assume that the returns are independent and identically distributed (i.i.d.) random variables satisfying
	$$
	X_{\min}  \leq X(k) \leq X_{\max}
	$$
	with $ -1 < X_{\min} < 0 < X_{\max} < \infty $ being known bounds.

	\vspace{5mm}
	\textbf{Frequency Considerations with Zero Execution Delay:}
	To study performance for a high-frequency trader, we use~$V_1(k)$ to denote the account value at stage $k$ and the investment takes the form of linear~feedback
	$$
	I_1(k) = K V_1(k)
	$$
	where~$K \in \mathcal{K} \doteq [0,1]$ is the fraction of the trader's account at risk.  
	For any admissible value of $K$, beginning with $V_1(0) = V(0)>0$,
	the dynamic evolution of the account value is described by the recursive~equation
	\begin{align*}
		V_1(k+1)    & = V_1(k) + I_1(k) X(k)\\[.5ex]
		& = (1 + K X(k) )V_1(k). 
	\end{align*}
	In the sequel, since we primarily focus on the final account value $V_1(n)$, whenever convenient, we use notation~$V_1(n,K)$ to emphasize the dependence on $K$ at this endpoint with~\mbox{$k=n$}.
	
	\vspace{3mm}
	On the other hand, for the buy and holder who does no rebalancing, using~$V_n(k)$ to denote the account value at stage $k$, with initial account value $V_n(0) = V(0)$, this trader uses~investment 
	$$
	I_n(0) = K V(0)
	$$
	where $K \in [0,1]$. 
	Note that the fraction $K$ used for the buy and holder is not necessarily the same $K$ used for the high-frequency trader.
	% by ordering
	%$$
	%N_0 \doteq \frac{K V(0)}{S(0)}
	%$$
	%shares. 
	The associated account value evolves over time~via
	\begin{align*}
		V_n(k+1)    &= V_n(k) + I_n(0)X(k) \\[.5ex]
		&= V_n(k) + K V(0)X(k).
	\end{align*}
	%$$
	%V(k+1) = V(k) + N_0(S(k+1)-S(k)).
	%$$
	% for $k > 1$ 
	In the sequel, similar to the case of the high-frequency trader, whenever convenient, we  use~$V_n(n,K)$ instead of~$V_n(n)$ to emphasize the dependence on $K$ at the endpoint~\mbox{$k=n$}. 
	% In addition, the associated solution to the recursive equation above is readily shown as
	%\begin{align*}
	%V(n) &= \left(1 + K\left( \prod_{k=0}^{n-1} (1+X(k)) -1\right) \right)V(0).  
	%\end{align*}

	\vspace{3mm}
	Now, for the two traders above,  we consider the expected logarithmic growths
	\begin{align*}
		g_1(K) &\doteq \frac{1}{n}\mathbb{E}\left[ \log \frac{V_1(n,K)}{V(0)}\right];\\[.5ex]
		g_n(K) &\doteq \frac{1}{n}\mathbb{E}\left[ \log \frac{V_n(n,K)}{V(0)}\right]
	\end{align*}
	%	where $g_1(K)$ corresponds to the performance achieved by high-frequency trading and  $g_{n}(K)$ corresponds to the performance achieved by using buy and hold.
	and our goal is  to find optima~$K_1^*$ and $K_n^*$ achieving
	\[
	g_1^* = \max_{K \in \mathcal{K}} g_1(K)
	\]
	and
	\[
	g_n^* = \max_{K \in \mathcal{K}} g_n(K)
	\] 
	respectively. In the sequel,  $K_1^*, K_n^* \in \mathcal{K}$ are called  \textit{optimal Kelly fractions} for high-frequency trading and buy and hold,~respectively.
	\vspace{8mm}
	\section{\small No Delay: Buy and Hold Versus High-Frequency}
	%
	% No Delay and the Maximality Theorem
	%
	\label{SEC: The Blessing Of High-Frequency Trading}
	\vspace{0mm}	
	In our previous paper~\cite{Hsieh_Barmish_Gubner_2018_ACC}, we introduced the notion of a  \textit{sufficiently attractive} stock; i.e., one whose i.i.d.\ returns~$X(k)$~satisfy
	$$
	\mathbb{E}\bigg[ \frac{1}{1+X(k)}\bigg]  \leq 1.
	$$ 
	We then proved that under this condition, the buy and holder matches the ELG performance achieved by high-frequency trader; i.e., $g_n^*=g_1^*.$
	%For example, if $X(k)$ takes the value~$X_{\max}$ with probability $p$ and $X_{\min}$ with probability~$1-p$, then the sufficient attractiveness condition holds if and only if
	%\[
	%p \geq \frac{X_{\min}(1+X_{\max})}{X_{\min}-X_{\max}}. 
	%\]
	%	We will revisit this condition later in Section~\ref{SEC: Example Include Execution Delay} when we provide our illustrative example. 
	
	\vspace{3mm}
	In the theorem below,
	%	 consistent with the formulation in the last section, to avoid trivialities, no transaction costs are assumed in the theorem below. 
	we prove that with no execution delay, high-frequency trading is \textit{unbeatable} in this same ELG content; i.e.,~\mbox{$g_n^* \leq g_1^*.$}  In Section~\ref{SEC: Example Include Execution Delay}, we see that this is no longer the case when execution delays and associated self-financing considerations are incorporated into the model.
	
	%	\vspace{3mm}
	%while it is straightforward to prove and intuitively obvious that the inclusion of transaction costs can negate the benefits of high-frequency trading, consistent with the formulation in the last section, to avoid trivialities, no transaction costs are assumed in the theorem below. 
	%We see that high-frequency trading is \textit{unbeatable} in terms of ELG; i.e., $g_n^* \geq g_1^*.$  

	%In this section, we prove  that with no execution delay and no transaction costs, high-frequency trading is unbeatable in terms of ELG; i.e., $g_n^* \geq g_1^*.$ 
	%
	% would point out that we had achieved g_n^* = g_1^* and wondered it g_n^* > g_1^* is possible. In the present paper we prove that it is not; i.e., we prove that g_n^* <= g_1^* when there is no delay and there are no transaction costs.
	
	%	\vspace{3mm}
	%	

	\vspace{5mm}
	\textbf{High-Frequency Maximality Theorem: }{\it 
		For the frequency-based trading scenario defined in Section~\ref{SEC: Problem Formulation}, it follows that
		\[
		g_n^* \leq g_1^*.
		\]
	}
	
	\vspace{0mm}
	\iftrue
	{\bf Proof:} Using shorthand $X_k$ for $X(k)$, for $K \in [0,1]$, the account value of the high-frequency trader is given by
	\[
	V_1(n,K) = \prod_{k=0}^{n-1} (1+K X_k)V(0).
	\]
	Note that since~\mbox{$X_k>-1$} for all~$k$, and since~\mbox{$0 \leq K \leq 1$}, we have~\mbox{$
		1+K X_{k}>0.
		$}
	Since the $X_k$ are i.i.d., the associated expected logarithmic growth is
	\begin{align*}
		g_1(K) &= \frac{1}{n}\mathbb{E}\left[ \log \frac{V_1(n,K)}{V(0)}\right]\\[.5ex] 
		&= \mathbb{E}[\log(1 +K X_{n-1})].
	\end{align*}

	The  account value for the buy and holder is given by
	$
	V_n(n,K) = \left(1 + K\mathcal{X}_n  \right)V(0)
	$,
	where 
	$$
	\mathcal{X}_n \doteq  \prod_{k=0}^{n-1} (1+X_k) -1
	$$
	is the compound return. To show~$g_n^* \leq g_1^*$, we use the smoothing property of conditional expectation to write the expected logarithmic growth as
	\begin{align*}
		g_n(K) &= \frac{1}{n}\mathbb{E}\left[ \log \frac{V_n(n,K)}{V(0)}\right]\\ &=\frac{1}{n} \mathbb{E}[\log(1+K\mathcal{X}_n)]\\[.5ex]
		&=\frac{1}{n} \mathbb{E}\left[\, \mathbb{E}[\log(1+K \mathcal{X}_{n})|X_{n-1}]\,\right]\\[.5ex]
		&= \frac{1}{n}\mathbb{E}\left[\, \mathbb{E}\left[\log\frac{1+K\mathcal{X}_{n}}{1+KX_{n-1}}\bigg|X_{n-1}\right]\,\right]\\[.5ex]
		& \hspace{5mm}+\frac{1}{n}\mathbb{E}\left[\, \mathbb{E}\left[\log(1+KX_{n-1})|X_{n-1}\right]\,\right]
		\\[.5ex]
		& =\frac{1}{n}\mathbb{E}\left[\, \mathbb{E}\left[\log\frac{1+K\mathcal{X}_{n}}{1+KX_{n-1}}\bigg|X_{n-1}\right]\,\right]\\[.5ex]
		&  \hspace{5mm}+\frac{1}{n} \mathbb{E}\left[\log(1+KX_{n-1})\right]
		\\[.5ex]
		& =\frac{1}{n}\mathbb{E}\left[\, \mathbb{E}\left[\log\frac{1+K\mathcal{X}_{n}}{1+KX_{n-1}}\bigg|X_{n-1}\right]\,\right]
		%\\[.5ex]
		%&  \hspace{5mm}
		+\frac{1}{n} g_1(K),
	\end{align*}
	where the last step follows from the calculation at the beginning of the proof.
	The simplification of the inner conditional expectation uses
	the independence of $\mathcal{X}_{n-1}$ and $X_{n-1}$ and the fact that
	$$
	\mathcal{X}_{n} - X_{n-1} =  (1+X_{n-1})\mathcal{X}_{n-1}.
	$$
	%$$
	%1+K\mathcal{X}_{n} = (1+KX_{n-1})+K\mathcal{X}_{n-1}(1+X_{n-1}),
	%$$
	Conditioning on $X_{n-1}=x$ for $x>-1$, we can write
	\begin{align*}
		&\mathbb{E}\left[\log\frac{1+K\mathcal{X}_{n}}{1+KX_{n-1}}\bigg|X_{n-1} = x\right]\\[.5ex]
		% &= \mathbb{E}\left[\log\left( 1+ \frac{K(1+X(n))}{1+KX(n)}\mathcal{X}_n \right)|X(n)\right]\\
		&\hspace{5mm} = \mathbb{E}\left[\log\left( 1+ \frac{K(\mathcal{X}_n - X_{n-1})}{1+KX_{n-1}}  \right)\bigg|X_{n-1}=x\right]\\[.5ex]
		&\hspace{5mm} = \mathbb{E}\left[\log\left( 1+ \frac{K(1+X_{n-1})}{1+KX_{n-1}}  \mathcal{X}_{n-1} \right)\bigg|X_{n-1}=x\right]\\[.5ex]
		&\hspace{5mm} = \mathbb{E}\left[\log\left( 1+ \frac{K(1+x)}{1+Kx}  \mathcal{X}_{n-1} \right)\bigg|X_{n-1}=x\right]\\[.5ex]
		&\hspace{5mm} = \mathbb{E}\left[\log\left( 1+ \frac{K(1+x)}{1+Kx}  \mathcal{X}_{n-1} \right)\right]\\[.5ex]
		&\hspace{5mm}= (n-1)g_{n-1}\left(K_x\right)
		%	&\hspace{5mm}\leq (n-1)g_{n-1}^*
	\end{align*}
	where
	$$
	K_x \doteq \frac{K(1+x)}{1+Kx}.
	$$
	Since $K \in [0,1]$ and $x>-1$, it follows that $K_x \in [0,1]$. Hence,
	$$
	g_{n-1}(K_x)\le \sup_{K\in[0,1]} g_{n-1}(K).
	$$
	Since the right-hand side is equal to $g_{n-1}^*$, it follows that
	\begin{align*}
		\mathbb{E}\left[\log\frac{1+K\mathcal{X}_{n}}{1+KX_{n-1}}\bigg|X_{n-1} = x\right] 
		%&= (n-1)g_{n-1}\left(K_x \right) \\
		&\leq (n-1)g_{n-1}^*.
	\end{align*}
	We now have
	\begin{align*}
		g_n(K) &\le \frac{n-1}{n}g_{n-1}^* + \frac{1}{n}g_1(K).
		%\\[.5ex]
		%& \le \frac{n-1}{n}g_{n-1}^* + \frac{1}{n}g_1^*.
	\end{align*}
	Taking the supremum over~\mbox{$K \in [0,1]$}, we obtain
	\[
	g_{n}^* \leq \frac{n-1}{n}g_{n-1}^* + \frac{1}{n}g_1^*.
	\]
	% We now observe that t
	To complete the proof, it is now noted that the foregoing argument for $g_n^*$ also applies to any $g_m^*$ for $m>1$. 
	Hence, 
	$$
	g_{m}^* \leq \frac{m-1}{m}g_{m-1}^* + \frac{1}{m}g_1^*.
	$$
	Now with $m=2$, we have $g_2^* \le g_1^*$. Similarly, for $m=3,$ it follows that
	$$g_3^* \le \frac{2}{3}g_{2}^* + \frac{1}{3}g_1^* \leq g_1^*,$$ 
	Continuing in this way we arrive at~\mbox{$g_n^* \le g_1^*$}. \hspace{5mm} $\square$	      
	\else
	{\bf Proof:} We begin by noting that the  account value for buy and hold is
	$
	V_n(n,K) = \left(1 + K\mathcal{X}_n  \right)V(0)
	$,
	where 
	$$
	\mathcal{X}_n \doteq  \prod_{k=0}^{n-1} (1+X(k)) -1
	$$
	is the compound return. In the sequel, to simplify the exposition, we use the shorthand notation $X_{n-1}$  for~$X(n-1)$.  To show~$g_n^* \leq g_1^*$, we first note that since~\mbox{$X(k)>-1$} for all~\mbox{$k = 0,1,\ldots, n-1$} and since~\mbox{$0 \leq K \leq 1$}, we have~\mbox{$
		1+K X_{n-1}>0.
		$}
	Next, using the smoothing property of conditional expectation, the expected logarithmic growth~becomes
	\begin{align*}
		g_n(K) &= \frac{1}{n}\mathbb{E}\left[ \log \frac{V_n(n,K)}{V(0)}\right]\\ &=\frac{1}{n} \mathbb{E}[\log(1+K\mathcal{X}_n)]\\[.5ex]
		&=\frac{1}{n} \mathbb{E}\left[\, \mathbb{E}[\log(1+K \mathcal{X}_{n})|X_{n-1}]\,\right]\\[.5ex]
		&= \frac{1}{n}\mathbb{E}\left[\, \mathbb{E}\left[\log\frac{1+K\mathcal{X}_{n}}{1+KX_{n-1}}\bigg|X_{n-1}\right]\,\right]\\[.5ex]
		& \hspace{5mm}+\frac{1}{n}\mathbb{E}\left[\, \mathbb{E}\left[\log(1+KX_{n-1})|X_{n-1}\right]\,\right]\\[.5ex]
		& =\frac{1}{n}\mathbb{E}\left[\, \mathbb{E}\left[\log\frac{1+K\mathcal{X}_{n}}{1+KX_{n-1}}\bigg|X_{n-1}\right]\,\right]\\[.5ex]
		&  \hspace{5mm}+\frac{1}{n} \mathbb{E}\left[\log(1+KX_{n-1})\right].
	\end{align*}
	Using the independence of $\mathcal{X}_{n-1}$ and $X_{n-1}$ and the fact that
	$$
	\mathcal{X}_{n} - X_{n-1} =  (1+X_{n-1})\mathcal{X}_{n-1},
	$$
	%$$
	%1+K\mathcal{X}_{n} = (1+KX_{n-1})+K\mathcal{X}_{n-1}(1+X_{n-1}),
	%$$
	we condition on $X_{n-1}=x$ for $x>-1$. Now  the first inner expectation above becomes 
	\begin{align*}
		&\mathbb{E}\left[\log\frac{1+K\mathcal{X}_{n}}{1+KX_{n-1}}\bigg|X_{n-1} = x\right]\\[.5ex]
		% &= \mathbb{E}\left[\log\left( 1+ \frac{K(1+X(n))}{1+KX(n)}\mathcal{X}_n \right)|X(n)\right]\\
		&\hspace{5mm} = \mathbb{E}\left[\log\left( 1+ \frac{K(\mathcal{X}_n - X_{n-1})}{1+KX_{n-1}}  \right)\bigg|X_{n-1}=x\right]\\[.5ex]
		&\hspace{5mm} = \mathbb{E}\left[\log\left( 1+ \frac{K(1+X_{n-1})}{1+KX_{n-1}}  \mathcal{X}_{n-1} \right)\bigg|X_{n-1}=x\right]\\[.5ex]
		&\hspace{5mm} = \mathbb{E}\left[\log\left( 1+ \frac{K(1+x)}{1+Kx}  \mathcal{X}_{n-1} \right)\bigg|X_{n-1}=x\right]\\[.5ex]
		&\hspace{5mm} = \mathbb{E}\left[\log\left( 1+ \frac{K(1+x)}{1+Kx}  \mathcal{X}_{n-1} \right)\right]\\[.5ex]
		&\hspace{5mm}= (n-1)g_{n-1}\left(K_x\right)
		%	&\hspace{5mm}\leq (n-1)g_{n-1}^*
	\end{align*}
	where
	$$
	K_x \doteq \frac{K(1+x)}{1+Kx}.
	$$
	Noting that $K_x \in [0,1]$ is  admissible, the equality above becomes
	\begin{align*}
		\mathbb{E}\left[\log\frac{1+K\mathcal{X}_{n}}{1+KX_{n-1}}\bigg|X_{n-1} = x\right] 
		%&= (n-1)g_{n-1}\left(K_x \right) \\
		&\leq (n-1)g_{n-1}^*.
	\end{align*}
	Next, 
	%using the fact that $X(k)$ are i.i.d., 
	we examine the second expectation term comprising~$g_n(K)$. Note that  via a straightforward calculation, the account value for high-frequency trading at stage $k=n$  is given by
	\[
	V_1(n,K) = \prod_{k=0}^{n-1} (1+KX(k))V(0).
	\]
	Hence, using the fact that $X(k)$ are i.i.d., the associated expected logarithmic growth is
	\begin{align*}
		g_1(K) &= \frac{1}{n}\mathbb{E}\left[ \log \frac{V_1(n,K)}{V(0)}\right]\\[.5ex] 
		&= \mathbb{E}[\log(1 +K X_{n-1})].
	\end{align*}
	Since $g_1(K) \leq g_1^*$,
	%	$
	%	\mathbb{E}[\log(1+KX_{n-1})] \leq g_1^*.
	%	$
	we arrive~at
	\begin{align*}
		%\mathbb{E}[\log(1+K \mathcal{X}_{n})]
		% &= \mathbb{E}\left[\, \mathbb{E}[\log(1+K \mathcal{X}_{m})|X_{m-1}]\,\right] \\[.5ex]
		%& =  \mathbb{E}\left[ \mathbb{E}\left[\log\frac{1+K\mathcal{X}_{m}}{1+KX_{m-1}}\bigg|X_{m-1}\right] \right] \\
		%& \hspace{5mm} + \mathbb{E}\left[ \mathbb{E}\left[\log\left( 1+ K X_{m-1} \right)\right| X_{m-1} \right] ]\\[.5ex]
		%& =  \mathbb{E}\left[ \mathbb{E}\left[\log\left( 1+ \widetilde{\kappa}\mathcal{X}_{n-1} \right)\bigg| X_{n-1}\right] \right]\\
		% & \hspace{5mm} + \mathbb{E}\left[\log\left( 1+ {K}{X}_{n-1} \right)\right]\\[.5ex]
		%& = (n-1)g_{n-1}(\widetilde{\kappa}) + g_1(K)\\[.5ex]
		g_n(K)&   
		%\frac{1}{n}\mathbb{E}\left[\, \mathbb{E}\left[\log\frac{1+K\mathcal{X}_{n}}{1+KX_{n-1}}\bigg|X_{n-1}\right]\,\right]\\[.5ex]
		%& \hspace{5mm}+\frac{1}{n}\mathbb{E}\left[\, \mathbb{E}\left[\log(1+KX_{n-1})|X_{n-1}\right]\,\right]\\[.5ex]
		\leq \frac{n-1}{n}g_{n-1}^*+\frac{1}{n}g_1^*.
	\end{align*}
	Taking the supremum over~\mbox{$K \in [0,1]$}, we obtain
	\[
	g_{n}^* \leq \frac{n-1}{n}g_{n-1}^* + \frac{1}{n}g_1^*.
	\]
	We now observe that the foregoing argument for $g_n^*$ applies to any $g_m^*$ for $m>1$. 
	Hence, 
	$$
	g_{m}^* \leq \frac{m-1}{m}g_{m-1}^* + \frac{1}{m}g_1^*.
	$$
	Now with $m=2$, we have $g_2^* \le g_1^*$. Similarly, for $m=3,$ it follows that
	$$g_3^* \le \frac{2}{3}g_{2}^* + \frac{1}{3}g_1^* \leq g_1^*,$$ 
	Continuing in this way we arrive at~\mbox{$g_n^* \le g_1^*$}. \hspace{5mm} $\square$
	\fi
	
	%\vspace{5mm}
	%{\bf Remarks:} It is worth mentioning that the inequality stated in the proof i.e.,
	%\[
	%g_{n}^* \leq \frac{n-1}{n}g_{n-1}^* + \frac{1}{n}g_1^*
	%\]
	%implies more than just $g_1^* \geq g_n^*$. In particular, it characterizes that the expected logarithmic growth $g_n^*$ has certain monotonicity in $n$.
	
	%\vspace{5mm}
	%\textbf{Sufficient Attractiveness Considerations:} 

	\vspace{8mm}
	\section{Extended Formulation with Execution Delay}
	\label{SEC: Extension to Include Delay In Trade Execution}
	Motivated by the fact that a trader's interactions with the market are not instantaneous, in this section, our aim  is to extend the formulation in Section~\ref{SEC: Problem Formulation} to incorporate a one-step delay in trade execution. One complication is that while in the delay-free case, an order specified in dollars is equivalent to an order specified in shares, when orders are delayed this is no longer true. This is an important observation because brokers typically accept orders in shares, not dollars.
	To illustrate the difference between orders in dollars versus shares, we first suppose that a broker were willing to accept an order in dollars.
	Then an order made at stage $k=0$ to buy shares worth~$K V(0)$ dollars is executed at stage~$k=1$ at price $S(1)$. This results in a stock holding at $k=1$, 
	%$K V(0)/S(1)$ 
	%shares, 
	whose value is exactly~$K V(0)$.
	In contrast, consider an order made at stage $k=0$ for
	$$
	N(0) \doteq \frac{K V(0)}{S(0)}
	$$
	shares, which is executed at stage $k=1$ at price $S(1)$. The value of these shares is
	%If there were no delay in executing the trade, the shares would be bought instantly at price $S(0)$, and their value would be $K V(0)$. However, the one-step delay in trade execution means that the $N(0)$ shares will be bought at stage~\mbox{$k=1$} at price $S(1)$, and their value will be
	$$
	N(0)S(1)=K V(0)(1+X(0)).
	$$
	Thus, depending on $X(0)$, orders in dollars versus shares can lead to very different results. In the sequel, following broker practices, orders are expressed in shares.
	It should be emphasized here that the analysis above at $k=0$ holds for both the high-frequency~trader and for the buy and~holder. 
	%In the sequel, we  work with shares for the order. 

	%	\vspace{5mm}
	%	\textbf{Remark:}
	%	If there were no delay, the amount at risk at stage~\mbox{$k=1$} would be
	%	$$
	%	\Bigl(\frac{K V_1(1)}{S(1)}\Bigr) S(1)
	%	= K V_1(1).
	%	$$
	%	This illustrates why in the no-delay case,
	%	it was not necessary to introduce shares.
	
	\vspace{3mm}
	For the case of the high-frequency trader, similarly, our convention is that at stage $k$, the trader
	places an order for
	$$
	N_1(k) \doteq \frac{K V_1(k)}{S(k)}
	$$
	shares	which are purchased at stage $k+1$ at price $S(k+1)$.
	%	Hence, the amount at risk is $N(k)S(k+1)$ at stage $k+1$.

	\vspace{5mm}
	{\bf Account Value Dynamics with Delay:}
	As in the zero-delay case described in Section~\ref{SEC: Problem Formulation},  the trade is required to be long only and self-financed. To be more specific, for the high-frequency trader, we require that the corresponding investment executed at stage $k$
	$$
	I_1(k) \doteq N_1(k-1)S(k) 
	%	=  KV_1(k-1)(1+X(k-1)),
	$$
	satisfies
	$
	0 \leq I_1(k) \leq V_1(k)
	$
	for all $k \geq 1$.
	Then the evolution of the account value  is described~by 
	\begin{align*}
		V_1(k+1) 
		&= V_1(k) + N_1(k-1)(S(k+1)-S(k))
		% V_1(2) = V_1(1) + N_1(0) (S(2) - S(1))  @k=1
		%        = ( 1 +  K X(1) S(1)/S(0) )V_1(0) 
		%        = ( 1 +  K X(1) (1+X(0)) )V_1(0) 
	\end{align*}
	for $k \geq 1$
	with~\mbox{$V_1(0)=V_1(1)= V(0)>0$}.  
	%	In the sequel, to emphasize the dependence on the fraction $K$, we use $V_{1}(n,K)$ instead of $V_1(n)$ for the solution of the recursive equation above at the endpoint~\mbox{$k=n$}. 

	\vspace{3mm}
	On the other hand, for the  buy and holder, since only  one order is executed at stage $k=1$, the 
	long-only and self-financing conditions force the corresponding investment,
	$$
	I_n(1) \doteq N_n(0)S(1)
	%	= KV_n(0)(1+X(0)).
	$$
	where $N_n(0) \doteq {K V_n (0)}/{S(0)}$,
	to satisfy
	$$
	0\leq I_n(1) \leq V_n(1).
	$$
	Then the corresponding account value is readily shown to satisfy the recursion
	$$
	V_n(k+1) = V_n(k) + N_n(0)(S(k+1)-S(k))
	$$
	for $k \geq 1$ with $V_n(0) = V_n(1) = V(0) >0$. Given the fact that the number of shares never changes, it is straightforward to obtain the closed-form
	\[
	V_n(n) = \biggl( 1 +K(1+X(0)) \biggl(\prod_{k=1}^{n-1} (1+X(k))-1\biggr) \biggr) V_n(0).
	\]
	%	Henceforth, similar to the high-frequency trading case, whenever convenient, we use $V_{n}(n,K)$ instead of $V_n(n)$ to emphasize the dependence on $K$ for the solution at the endpoint~\mbox{$k=n$}. 
	Similar to the case without delay, for ELG purposes, we  use the notation~\mbox{$g_{1}(K)$} and $g_{n}(K)$ to denote the performance, as a function of~$K$, achieved by high-frequency trading and buy and hold, respectively.
	In addition, we denote optima by~$K_1^*$ and $K_n^*$ and the associated optimal values by~$g_1^*$ and~$g_n^*$.
	
	%%%%%%%%%%%%%%%%%%%%%%%%%%%
	\vspace{5mm}
	\textbf{On Long-Only and Self-Financing with Delay:} When execution delay is in play, in contrast to the no-delay case,~\mbox{$K \leq 1$} does not guarantee self-financing. 
	%	As the no-delay case, this requirement is equivalent to characterize the range of~$K$.
	For the buy and holder at stage $k=1$, self-financing requires~\mbox{$I_n(1) \leq V_n(1),$} 
	which is equivalent~to
	\[
	\frac{K V_n (0)}{S(0)}S(1) \leq V_n(1) .
	\]	
	Now using the fact that~\mbox{
		$
		{S(1)}/{S(0)} = 1+X(0)
		$} and $V_n(0)=V_n(1)$,
	%, we obtain~\mbox{$0\leq   {K} (1+X(0)) \leq 1 $}. To assure the inequality holds with probability one, it is equivalent to require
	the inequality above holds for all possible values of $X(0)$ if and only if 
	$$
	K    \leq \frac{1}{1+X_{\max}}.
	$$
	Combining this with the long only  constraint that $I_n(1) \geq 0$,
	we have
	$$
	0 \le K \le \frac{1}{1+X_{\max}}.
	$$
	For the case of the high-frequency trader, as seen in the lemma below, once again, the same restriction on $K$ results, but a lengthier argument is required.
	%It is important to note that if the inequality above on~$K$  holds, a  induction argument leads to the fact that  bankruptcy never occurs for both the high-frequency trader and the buy and holder; i.e.,~$V_i(k) \geq 0$ for all $k \geq 1$ and~\mbox{$i \in \{1,n\}$.} We use this fact implicitly in the following lemma.

	\vspace{5mm}
	{\bf The Self-Financing Lemma: }{\it For the case of one-step delay in execution,  the high-frequency trader is long only and self-financed if and only if
		$$
		0\leq K \leq \frac{1}{1+X_{\max}}. 
		$$ 
		Furthermore, when $K$ satisfies the above inequality, the trader's account is nonnegative; i.e., $V_1(k) \geq 0$ for all~\mbox{$k\geq 0.$}
		%		$$
		%		0\leq I_1(k) \leq V_1(k)
		%		$$
		%		for all $k \geq 1$  if and only if the inequality on $K$ above~holds. 
	}

	\vspace{5mm}
	{\bf Proof:}
	%We begin by showing that the bounds on $K$ are necessary and sufficient for the investment at $k=1$ to be long only and self-financed That is, $0 \leq I_1(k) \leq V_1(k)$. 
	%Then we show that the bounds on $K$ are sufficient to guarantee that the investments for $k\geq 1$ are long only and self-financed.
	%\vspace{3mm}
	The necessity of the conditions on $K$ follows by the same argument used for the buy and holder given preceding the lemma statement.
	\vspace{3mm}
	%Since 
	%\[
	%I_1(1) 	= K(1+X(0))V_1(0) 
	%\]
	%and since $1+X(0)>0$ and $V_1(0)>0$, we see that~$I_1(1) \geq 0$ if and only if $K \geq 0$. Since $V_1(0) = V_1(1)$, we see that~$I_1(1) \leq V_1(1)$ if and only if $K \leq 1/(1+X_{\max})$.  In particular, we have shown the necessity of the bounds on~$K$.
	%This also implies that the proof of necessity is completed.
	
	%	To prove necessity, we proceed a proof by contradiction. Specifically, we show that if 
	%	% assume~\mbox{$I_1(k) \leq V_1(k)$} holds along all sample paths, we must show the required condition on $K$ holds. Proceeds by contradiction, suppose 
	%	$
	%	K > {1}/({1+X_{\max}})
	%	$
	%	then $I_1(1)>V_1(1)$. Indeed, when $X(0)=X_{\max}$ and~\mbox{$K>1/(1+X_{\max})$}, we have 
	%	\begin{align*}
	%	I_1(1) 	&= K(1+X(0))V_1(0) 
	%	>   V_1(0) = V_1(1). 
	%	\end{align*}
	%	This completes the proof of necessity. 
	
	To prove sufficiency, we assume that $K$ satisfies the given inequality. We also assume that the  $V_1(k)\ge 0$ noting that the argument that this holds is given at the end of the proof.
	%	then $V_1(k) \ge 0$ for $k\ge 1$, which can be proved by induction.
	Next, since
	\begin{align*}
		I_1(k) 	&\doteq  N_1(k-1){S(k)}\\[.5ex]
		&=(1+X(k-1))K V_1(k-1)
	\end{align*}
	must be nonnegative,
	% Assuming that the inequality on~$K$ holds, 
	it remains to show that $I_1(k) \leq V_1(k)$ for all~\mbox{$k \geq 1$}.
	Proceeding by induction, we begin by noting that for~$k=1$, 
	%	using the assumed bound on $K$,  the initial investment satisfies
	\begin{align*}
		I_1(1)  &= (1+X(0))K V_1(0) \\[.5ex]
		&\leq   (1+X_{\max})K V_1(0)  \leq V_1(0) = V_1(1).
	\end{align*}
	We next fix any $k \geq 1$, and suppose that for all paths~\mbox{$(X(0),X(1),\ldots,X(k-1) ),$}
	we have
	$
	I_1(i) \leq V_1(i)
	$ for~\mbox{$i \leq k$}.
	We must show $
	I_1(k+1) \leq V_1(k+1)
	$.
	We split the proof into two cases: 
	\\
	\textit{Case 1:}
	If $X(k)\geq 0$, then using the assumed bound on $K$ and the fact that  $I_1(k) \geq 0$, we obtain
	\begin{align*}
		I_1(k+1) 
		%&= (1+X(k))KV(k) \\
		&\leq (1+X_{\max})K V_1(k)\\[.5ex]
		& \leq V_1(k)\\[.5ex]
		&\leq  V_1(k) + I_1(k)X(k) = V_1(k+1).
	\end{align*}
	\textit{Case 2:} If $X(k)<0$, then, with the aid of the assumed inductive hypothesis $I_1(k)\leq V_1(k)$, we have
	\[
	I_1(k)X(k) \geq V_1(k)X(k).
	\]
%	\begin{align*}
%		I_1(k+1) &= (1+X(k))K V_1(k) \\[.5ex]
%		%	& \leq K V_1(k)\\[.5ex]
%		%	&\leq \frac{1}{1+X_{\max}} V_1(k)\\[.5ex]
%		& \leq V_1(k).
%	\end{align*}
Now, using the facts that $0\leq K \leq 1/(1+X_{\max}) < 1$ and $X(k)>-1$, we  observe that
	\begin{align*}
		V_1(k+1) &= V_1(k) + I_1(k)X(k) \\[.5ex]
			& \geq V_1(k) + V_1(k)X(k)\\[.5ex]
			& = 1\cdot (1+X(k))V_1(k)\\[.5ex]
			& > K(1+X(k))V_1(k)\\[.5ex]
			& = I_1(k+1).
	\end{align*}
% 
%	To show that~\mbox{$I_1(k+1) \leq V_1(k+1)$}, we
%	proceed by contradiction. If~\mbox{$I_1(k+1)>V_1(k+1)$},  it follows that~\mbox{$V_1(k) > V_1(k+1)$}. However, $X(k)<0$ in combination with $I_1(k) \geq 0$ assures that $V_1(k+1)<V_1(k)$
%	%$V_1(k+1) = V_1(k)+I_1(k)X(k) < V_1(k)$
%	which is the contradiction. 
	 This completes the proof of sufficiency.
	%	Therefore,~\mbox{$I_1(k+1) \leq V_1(k+1)$} is established.
	
	\vspace{3mm}
	To complete the proof of the lemma, it remains to show $V_1(k) \geq 0$ for all $k$. Noting that $V_1(0) = V_1(1) >0$, using  the assumed inequality on  $K$, we first see that 
	\begin{align*}
		V_1(2)  &= V_1(1) + N_1(0) (S(2) - S(1)) \\[.5ex]
		%	&= V_1(1) + I_1(1)X(1)\\
		&= ( 1+ K(1+X(0))X(1))V_1(0) \\[.5ex]
		& \geq ( 1+ \frac{1}{1+X_{\max}}(1+X_{\max})X_{\min})V_1(0)\\[.5ex] 
		&=  ( 1+ X_{\min})V_1(0) \geq 0.
	\end{align*}
	Then, continuing with an induction argument similar in flavor to the one above, it follows that $V_1(k) \geq 0$ for all~$k$.
	\hspace{5mm} $\square$

	%	\vspace{5mm}	
	%	\textbf{Remarks:}
	%	In the sequel, we use $\mathcal{K}_{SF}$ to denote
	%It is also worth mentioning that  if the required inequality above on $K$  holds, then bankruptcy never occurs for high-frequency trader; i.e.,~$V_1(k) \geq 0$ for all $k \geq 1$ with probability one.
	%  See the State Positivity Lemma for a detailed discussion. 
	
	%%%%%%%%%%%%%%%%%%%%%%%%%%%
	\vspace{8mm}
	\section{Delay: Buy and Hold Versus High-Frequency}
	\label{SEC: Example Include Execution Delay}
	\vspace{0mm}
	In this section, we show that trade execution delay can lead to better performance for a buy and holder versus that of the high-frequency trader.
	To accomplish this, we provide examples involving a binary lattice model for the stock returns. For such a model,~$X(k)$ takes the value $X_{\max}$ with probability~$p$ and the value~$X_{\min}$ with probability $1-p$.
	The rationale for use of the binary lattice is that the computations to follow are not too complex and that this model is used in finance. In addition this model also has the property that as the time $\Delta t$ between stages becomes small, one obtains an approximation of  classical Geometric Brownian Motion which is  widely used in the financial community;
	%, can be viewed as the limiting lattice behavior with~$\Delta t$ very small;
	e.g., see~\cite{Luenberger_2013}.
	
	\vspace{3mm}
	
	Our theoretical results in the no-delay case apply to a wide class of distributions of the returns. In particular, distributions that approximate real market returns are covered. In contrast, in this section on delay, we work with specific examples to make inferences about real markets. Hence, it is important that the distributions we use approximate what is seen in~practice.
	%In particular, for stock prices sampled about one second apart, maximal returns are on the order of $\pm 0.001$.

	\vspace{3mm}
	Before we provide our main example with $n=100$ steps and with returns that are a somewhat reasonable facsimile of real-world trading, we first analyze a toy example with only three trades and unrealistic returns. For this simple case, it is easy to show mathematically, rather than by simulation, that execution delay in combination with the self-financing requirement leads to 
	$$
	g_n^* >g_1^*.
	$$ 
	To this end, we let~\mbox{$n=3$} and use returns~\mbox{$X_{\max}=0.8$} and~\mbox{$X_{\min}=-0.2$} with equal probability. 
	%Note that,  when there is no execution delay,  per discussion in Section~\ref{SEC: The Blessing Of High-Frequency Trading}, it is readily verified that the returns generated from the binary model satisfy the ``sufficient attractiveness" inequality,
	%$
	%\mathbb{E}[ {1}/({1+X(k)})] \approx 0.9028 \leq 1.
	%$
	%Thus,  without execution delay, according to the Sufficiency Theorem in \cite{Hsieh_Barmish_Gubner_2018_ACC},  
	%Hence, the optimal fractions for both the high-frequency trader and the
	%buy and holder are~\mbox{$K_1^* = K_n^* = 1$} which leads to an identical optimal ELG; i.e.,~$g_1^* = g_n^*$.
	Since $n$ is small, a straightforward calculation allows one to obtain both~$g_1(K)$ and~$g_n(K)$ in closed form. First restricting $K$ to guarantee self-financing,
	%constraint imposed, according to the Self-Financing Lemma in Section~\ref{SEC: Extension to Include Delay In Trade Execution},
	we find that 
	$$
	K_1^* = K_n^* = \frac{1}{1+X_{\max}} \approx 0.556
	$$
	with associated ELGs given by are $g_1^* \approx 0.1009$ and~\mbox{$g_n^* \approx 0.1104$}. Hence, the buy and holder outperforms the high-frequency trader by about $9.44\%$. If the self-financing constraint is removed, say by allowing \mbox{$K \in [0, 1]$}, then a straightforward calculation leads to optimal fractions~\mbox{$K_1^* = K_n^* = 1$}, which corresponds to allowing \textit{leverage}. That is, the associated optimal investments satisfy
	\begin{align*}
		I_i^*(k) &= K_i^*(1+X(k))V_i(k) \\[.5ex]
		&\leq (1+X_{\max})V_i(k) = 1.8 V_i(k)
	\end{align*} 
	for $i\in \{1,n\}$.
	In this case, the associated ELGs are~\mbox{$g_1^* \approx 0.1237$} and $g_n^* \approx 0.1262$. Hence, the buy and holder outperforms the high-frequency trader by about~$2.11\%$. 
	Therefore, in this example, when delay is present, we see that if one drops the self-financing constraint, the buy and holder still outperforms the high-frequency trader. 
	This shows that delay alone, rather than in combination with the self-financing constraint, can allow the buy and holder to outperform the high-frequency~trader.

	\vspace{3mm}
	% It is important to note that lattice models are relevant to real stock market because the classical Geometric Brownian Motion, a widely used stock pries model, is obtained as its limiting behavior with~$\Delta t$ very small; e.g., see~\cite{Luenberger_2013}. 
	Since our goal is to argue that when execution delay is present, real markets might also see that the buy and holder outperforms the high-frequency trader, in the main example below, we work with $n=100$ rather than $n=3$, and smaller returns are used.
	Note that although the returns used below are small, they are not quite as small as those in real markets; see the discussion on this issue in Section~\ref{SEC: Conclusion and Future Work}.

	%	\vspace{0mm}
	%	\begin{center}
	%		\graphicspath{{figs/}}
	%		\includegraphics[scale=0.40]{latency_d=1_n=3_betting_v03.eps}
	%		\figcaption{Expected Log-Growth with $n=3$ and Large Returns}
	%		\label{fig:Expected Logarithmic Growth with $n=3$_binary}
	%	\end{center}
	%	\vspace{5mm}		

	%
	\vspace{5mm}
	{\bf Example (Binary Lattice Model):}  We consider the binary lattice model with returns $X_{\max} =  0.02$ with probability~\mbox{$p=0.6$} and~$X_{\min} = -0.01$ with probability~\mbox{$1-p=0.4$.}
	When there is no delay, we recall the sufficient attractiveness inequality from Section~\ref{SEC: The Blessing Of High-Frequency Trading} and note that
	for the more general binary lattice model parameterized in $X_{\min}$, $X_{\max}$ and $p$, the inequality reduces to
	$$
	p \ge \frac{X_{\min}(1+X_{\max})}{X_{\min}-X_{\max}}.
	$$
	For the lattice with $X_{\min}=-0.01$ and $X_{\max}=0.02$ under consideration, the sufficient attractiveness condition reduces to the requirement that $p \geq 0.34.$ Since the assumed value is $p=0.6$, the requirement is therefore satisfied.
	% 	which holds since the lattice under consideration has~\mbox{$p=0.6$} and the right-hand side is $0.34$. 
	Starting with~\mbox{$V(0)=10,000$} and stopping at stage~\mbox{$n = 100$},
	%when there is no execution delay,  per discussion in Section~\ref{SEC: The Blessing Of High-Frequency Trading}, it is readily verified that the returns   satisfy the ``sufficient attractiveness" inequality; i.e.,  , the inequality reduces
	%	$$
	%	\mathbb{E}\bigg[ \frac{1}{1+X(k)}\bigg] \approx 0.99 < 1.
	%	$$ 
	% reduces to
	%\[
	%p > \frac{X_{\min}(1+X_{\max})}{X_{\min}-X_{\max}} = 0.34.
	%\]
	we obtain the optimal fractions~$K_1^* = K_n^* = 1$ and identical optimal expected logarithmic growths; i.e.,~\mbox{$g_1^* = g_n^*$.}
	
	\vspace{3mm}
	For this same stock model example, we now consider the effect of a unit execution delay. According to the Self-Financing Lemma in Section~\ref{SEC: Extension to Include Delay In Trade Execution}, we require~\mbox{$
		0\leq K \leq {1}/{(1+X_{\max})} \approx 0.9804.
		$}
	Now, for the buy and holder, we use 
	the closed-form solution 
	in Section~\ref{SEC: Extension to Include Delay In Trade Execution} to calculate~$g_n(K)$. Indeed, a lengthy but straightforward calculation leads to
	%	We claim that
	%That is, use shorthand $Z_n \doteq \prod_{k=1}^{n-1} (1+X(k))-1$ and $X_0 \doteq X(0)$, we have
	\begin{align*}
		g_n(K)
		%		 &= \frac{1}{n} \mathbb{E} \bigg[ \log \frac{V_n(n,K)}{V(0)} \bigg] \\[.5ex]
		%&=\frac{1}{n} \mathbb{E} \bigg[ \log \left( 1 +K(1+X_0 ) Z_n \right) \bigg]\\[.5ex]
		%&=\frac{1}{n} p\mathbb{E} \bigg[ \log \left( 1 +K(1+X_0 ) Z_n\right) \bigg| X_0 =X_{\max} \bigg] \\
		%&\hspace{5mm} +\frac{1}{n}(1-p)\mathbb{E} \bigg[ \log \left( 1 +K(1+X_0 ) Z_n\right) \bigg| X_0 =X_{\min}\bigg]\\[.5ex]
		&=\frac{1}{n} p\sum_{i=0}^{n-1} p_i \log \left( 1 +K(1+X_{\max} )z_i \right) \\
		&\hspace{8mm} +\frac{1}{n}(1-p)\sum_{i=0}^{n-1} p_i \log \left( 1 +K(1+ X_{\min} )z_i \right)
	\end{align*}
	where
	$
	z_i \doteq (1 + X_{\max})^i(1 + X_{\min})^{n-1-i} - 1
	$
	and
	$$
	p_i \doteq \begin{pmatrix} n-1 \\ i \end{pmatrix} p^i (1-p)^{n-1-i}.
	$$ 
	%	with $X_{\min} = -0.01$, $X_{\max} = 0.02$ and $p=0.6$.
	%	This can be seen by using the law of total
	%	probability if we~define 
	%	$$
	%	Z_n \doteq \prod_{k=1}^{n-1} (1+X(k))-1
	%	$$
	%	and note that $X(k)$ takes only the
	%	values $X_{\min}$ and $X_{\max}$ and that $P(Z_n=z_i)=p_i$.
	Then, by plotting $g_n(K)$ versus $K$, we see in Figure~\ref{fig:Expected Logarithmic Growth with $n=100$_binary} that the optimal fraction~\mbox{$K_n^* \approx 0.9804$} corresponds to the limit imposed by self-financing. We also obtain the associated optimal ELG $g_n^* \approx 0.007719$; see the dash-dotted line in~Figure~\ref{fig:Expected Logarithmic Growth with $n=100$_binary}.

	\vspace{3mm}
	On the other hand, for the high-frequency trader, since a closed-form for $g_1(K)$ is unavailable, we perform a Monte-Carlo simulation using $500,000$ sample paths. In Figure~\ref{fig:Expected Logarithmic Growth with $n=100$_binary}, from the plots of~$g_1(K)$ and $g_n(K)$ versus~$K$, we obtain~\mbox{$K_1^* = K_n^* \approx 0.9804,$} which leads to the optimal expected logarithmic growth~\mbox{$g_1^* \approx 0.0076$.} Recalling that $g_n^* \approx 0.007719$, the optimal ELG for the buy-and-hold strategy exceeds that of the high-frequency trading strategy by about~$1.1\%$. The difference  $g_n^* - g_1^* > 0$ is consistently observed when one carries out many repetitions of the simulation. It is also noted that, if one drops the self-financing constraint, then the optimal fractions become~\mbox{$K_1^* = K_n^* = 1$}, which corresponds to allowing leverage as we saw in the~\mbox{$n=3$} case. That is, the optimal investments satisfy
	\begin{align*}
		I_i^*(k) &= K_i^*(1+X(k))V_i(k) \\[.5ex]
		&\leq (1+X_{\max})V_i(k) = 1.02 V_i(k)
	\end{align*} 
	for $i\in \{1,n\}$.
	In this case, we obtain~\mbox{$g_1^* \approx 0.0077$} and $g_n^* \approx 0.007826 $, which shows that the buy-and-hold strategy outperforms the high-frequency strategy by about~$0.56\%$.  This shows again that delay alone, rather than in combination with the self-financing constraint can lead to~\mbox{$g_n^*>g_1^*$.}
	
	\vspace{0mm}
	\begin{center}
		\graphicspath{{figs/}}
		\includegraphics[scale=0.40]{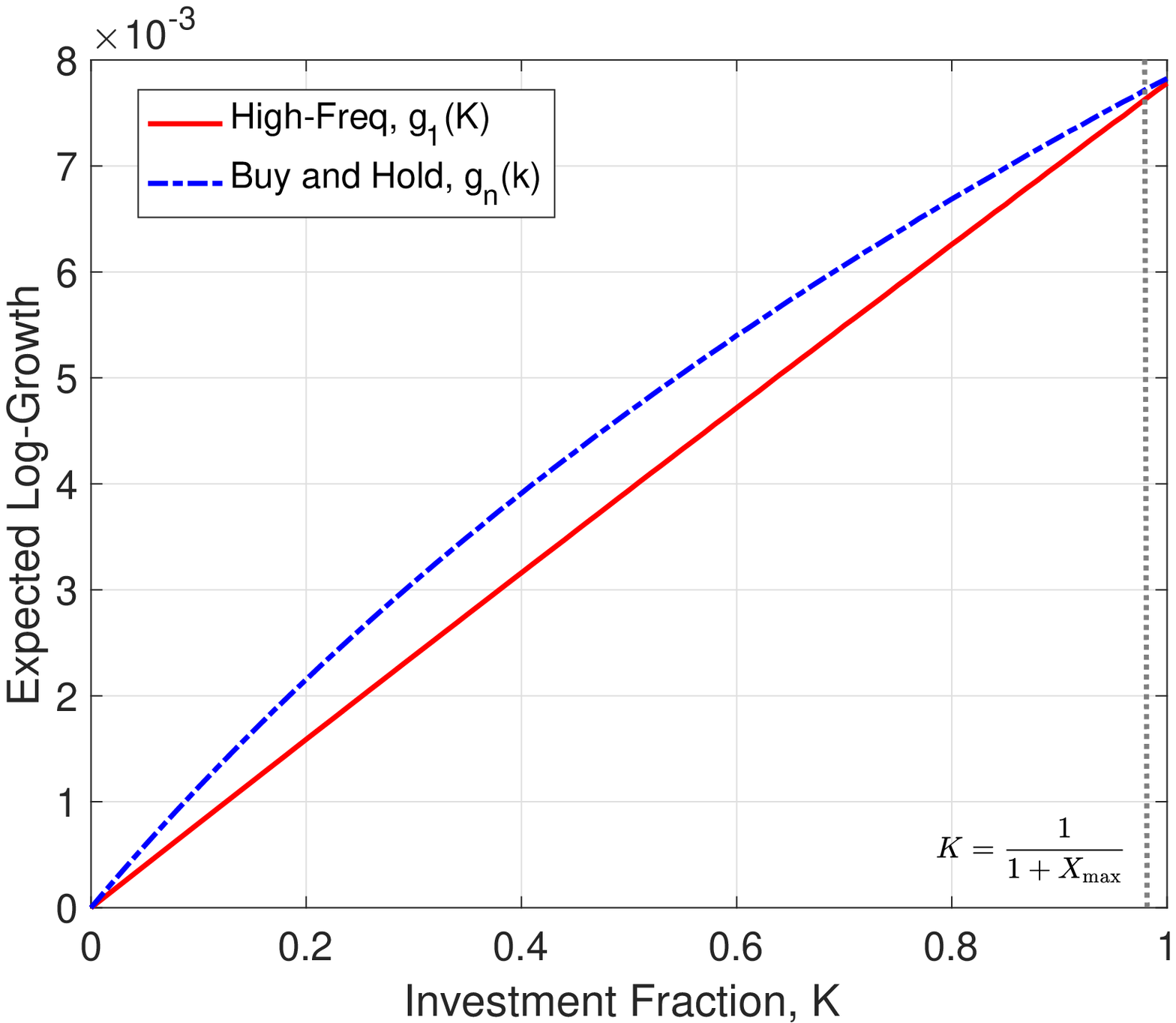}
		\figcaption{Expected Log-Growth with $n=100$}
		\label{fig:Expected Logarithmic Growth with $n=100$_binary}
	\end{center}
	\vspace{-2mm}

	%	\vspace{3mm}
	%	 \textbf{What If One Removed Self-Financing Constraint?} Moreover, even if self-financing condition is removed in the unit execution delay formulation, say one allows $K \in [0,1]$ rather than $K \in [0, 1/(1+X_{\max})]$,  then from Figure~\ref{fig:Expected Logarithmic Growth with $n=100$_binary}, we see optimal fractions for both high-frequency trader and buy and holder becomes one and buy and holder still has a margin of victory. Specifically, the ELGs we found are ~$g_1^* = g_1(1) \approx 0.0077$ and $g_n^* = g_n(1) \approx 0.007826 $.  Hence, the buy-and-hold strategy exceeds that of the high-frequency strategy by about~$0.56\%$.

	\vspace{5mm}
	{\bf Fractional Kelly Strategies:} The optimum above for~\mbox{$p=0.6$} requires that almost all funds be invested in the underlying stock. Since this might be viewed as  far too aggressive for many traders, many authors suggest using
	a so-called \textit{fractional} Kelly strategy. This is obtained  by scaling down the  fraction~$K$ so that the investment level is lower; e.g., see~\cite{Thorp_2006}, \cite{Maclean_Thorp_Ziemba_2010},  \cite{Rujeerapaiboon_Barmish_Kuhn_2018} and \cite{Davis_Lleo_2010}.  Now if one uses a fractional Kelly strategy for the binary lattice model above, as seen in Figure~\ref{fig:Expected Logarithmic Growth with $n=100$_binary}, the ``margin of victory" for the buy and holder can be larger. In fact,~\mbox{$g_n(K) > g_1(K) $} for the entire open interval $0<K<1$.
	%over the entire range of $K$ which is allowed by the self-financing~constraint. 
	% 	for any~\mbox{$0< K < 1/(1+X_{\max})$}. 

	\vspace{5mm}
	{\bf Binary Lattice with Variable Probability:} We now revisit the binary lattice example above with the same parameters~\mbox{$X_{\max}=0.02$}, $X_{\min} = -0.01$ and $n=100$, but now let the probability  $p$ vary. 
	Recalling the analysis in the previous subsection, when there is no delay,  the trade is sufficiently attractive if 
	%	\[
	%	\mathbb{E}\left[\frac{1}{1+X(k)}\right] \leq 1,
	%	\]
	%	which, after lengthy but straightforward calculation, is equivalent to
	$
	p >
	%	 \frac{X_{\min}(1+X_{\max})}{X_{\min}-X_{\max}} = 
	0.34.
	$
	Thus, within this range of $p$, for the no-delay case, the  buy and holder matches the performance of the high-frequency trader. However, when a unit execution delay is in play,
	%	Figure~\ref{fig:g_star_vs_p}, a plot of the  difference \mbox{$g_n^* - g_1^*$}  versus~\mbox{$p \in [0,1]$} indicates that  when~\mbox{$p \ge 0.34$}, we have~\mbox{$g_n^* - g_1^* >0$}. That is, the the buy and holder strictly outperforms the high-frequency~trader. Moreover, 
	Figure~\ref{fig:g_star_vs_p_in_percent} shows that the buy and holder becomes the ``winner." The plot of the percentage  difference \mbox{$(g_n^* - g_1^*)/g_1^* \times 100 \%$}  versus~\mbox{$p$} in this figure shows an increasing ``margin of victory" for the buy and holder as $p$ varies over its range.
	
	\vspace{0mm}
	\begin{center}
		\graphicspath{{figs/}}
		\includegraphics[scale=0.40]{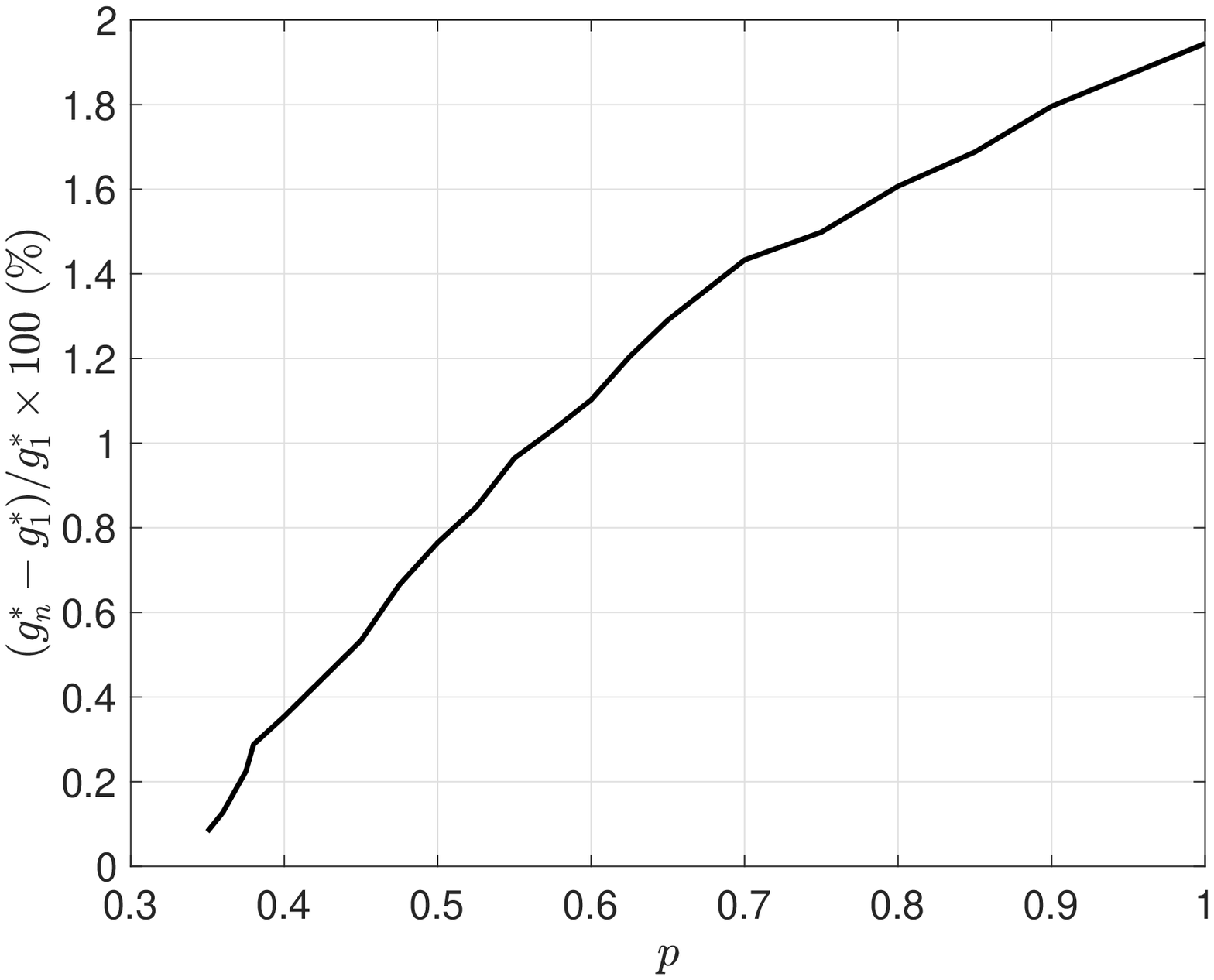}
		\figcaption{$(g_n^*-g_1^*)/g_1^* \times 100\%$ Versus $p$}
		\label{fig:g_star_vs_p_in_percent}
	\end{center}

	%	\vspace{5mm}
	%	\begin{center}
	%		\graphicspath{{figs/}}
	%		\includegraphics[scale=0.40]{binary_g_star_vs_p_diff.eps}
	%		\figcaption{$g_n^*-g_1^*$ Versus $p$ \\  }
	%		\label{fig:g_star_vs_p}
	%	\end{center}

	%
	%
	%
	%
	%\textbf{Remarks:} with large returns and small n leading to ambiguous victory for BH even when self-financing is not in play. 

	\vspace{8mm}
	\section{Conclusion and Future Work}
	\label{SEC: Conclusion and Future Work}
	In this paper, we studied a stock trading problem using Kelly's expected logarithmic growth criterion as the performance metric. 
	We first proved a theorem showing that high-frequency trading is unbeatable  when there are no transaction costs and no execution delays.
	Then, when delay in execution is considered, we showed, using a binary lattice model, that  there are cases when trading at high-frequency can be  inferior to buy and hold. 
	% For real market with returns really small, it is an open question whether buy and holder being the winner is realizable or not. 
	The binary lattice stock model used in Section~\ref{SEC: Example Include Execution Delay} is based on a mathematical model with returns~$X(k)$ which are larger than those seen with real ``high-frequency" trading data. Thus, it is important to ask whether the buy and holder can still achieve a higher ELG than a high-frequency trader based on a real-world model obtained from historical data. This situation is discussed~below.  
	
	\vspace{5mm}
	\textbf{Testing with Historical Data:} In this subsection, we provide a gateway to future research by describing how the ideas in this paper might be further pursued using historical intra-day tick data instead of a mathematical price model. For such data, each ``tick" corresponds to a new stock price, and the time between order book transactions can be as low as a microsecond. Let $t$ and $s$ denote the arrays of time stamps and transaction prices, respectively. For a trader confronted with an execution delay of $\Delta t$, in seconds the prices used for simulation are selected as follows:
	If a trade occurs at~$t(k)$, the next trading time which one should use occurs at time~$t(k')$~where 
	$$
	k' = \min\{i: t(i) \geq t(k) + \Delta t \}.
	$$
	With the convention above, we obtain a subsequence $t(k_j)$ of~$t(k)$ with total number of trades denoted by $m$,  and we
	%	Letting  $k_j$ with $j=0, \ldots ,m-1$ be the index for the subsequence of $t(k)$ described above.
	work with the \textit{associated} stock prices~$s( k_j )$ and the corresponding returns
	$$
	x_j \doteq \frac{s( k_{j+1} ) - s( k_{j} )}{s(k_{j} )}.
	$$
	These returns define an empirical probability mass function given by the sum of impulses
	$$
	\widehat{f}_{_X}(x) = \frac{1}{m} \sum_{j=0}^{m-1} \delta(x - x_j),
	$$
	which is used to generate random returns for the Monte-Carlo simulation which is needed to maximize expected logarithmic growth going~forward. 
	%	optimal Kelly fractions
	%	$K_1^*$, $K_n^*$, and associated optimal
	%	values~\mbox{$g_1^*$} and~\mbox{$g_n^* $.} 

	\vspace{3mm}
	Using the procedure above with $\Delta t = 1$, we conducted a preliminary experiment using  high-frequency historical intra-day tick data for APPLE \mbox{(ticker AAPL)}  for the period 9:30:00~AM to 10:00:00~AM on December~2,~2015.
	% 	Enforcing a one-second minimum time between executions~($\Delta t=1$), w
	We obtained $m=1293$ trades with approximately~$1.4$ seconds as the mean value of $t(k_{j+1})-t(k_j)$.
	Our Monte-Carlo simulation, carried out using $50,000$ sample paths for each value of $K$, resulted in optimal fractions
	are~\mbox{$K_1^* = K_n^* \approx 0.9978$}
	and associated optimal ELGs 
	\mbox{$g_1^* \approx 3.9966 \times 10^{-6}$} and~\mbox{$g_n^* \approx 4.002\times 10^{-6}$}. Since the difference between these quantities is too small to rule out roundoff error, we deemed our simulation to be inconclusive as to which of the two traders achieves better performance. 
	%	Although~$g_n^*>g_1^*$, the difference is very small, and it may be
	%	due to roundoff error or simulation variability.
	Another confounding factor to mention is that in our simulation,
	both traders own between
	$85.1$ and~$85.5$ shares
	%	; i.e., their holdings differ by small fractions of a share. 
	%	Based on these observations, we believe the conclusion that in practice, execution delay will often result in suboptimal performance of high-frequency trading compared with buy and hold cannot be made definitively at this point.	
	during the entirety of the thirty minutes of market time under consideration. In conclusion, studies of ELG performance using real market data are relegated to future research.
	
	\vspace{5mm}
	\textbf{Other Future Research Directions:}
	Regarding further research on delay-related issues, 
	%one obvious direction would be to see if the phenomenon that buy and holder can outperform high-frequency trader  with real world high-frequency trading intra-day tick data.
	another  interesting direction would be to extend our results  in trade execution to also include delay in information acquisition. 
	A second additional direction of future research is motivated by the fact that we only dealt with a single stock in this paper. In the future, we envision a formulation which involves a portfolio with multiple stocks. Our preliminary work to date suggests that a generalization of the High-Frequency Maximality Theorem given in Section~\ref{SEC: The Blessing Of High-Frequency Trading} should be possible to obtain. If this proves to be true, such a result would serve as a good starting point for future work.
	%	 Unlike in this paper where we only worked with a single stock in single execution delay case, we envision a formulation involves multiple delay is possible to establish. 
	%	 The third direction for future research involves extension of some of the results in this paper to a portfolio scenario; e.g., we envision that the High-Frequency Maximality Theorem in Section~\ref{SEC: The Blessing Of High-Frequency Trading} can be immediately extended to portfolio case.

	%%%%%%%%%%%%%%%%%%%%%%%%%%%%%%%%%%%%%%%%%%%%%%%%%%%%%%%%%%%%%%%%%%%%%%%%%%%%%%%%


\vspace{8mm}
	\begin{thebibliography}{99}
		\vspace{0mm}
		
		
		\bibitem{Kelly_1956}
		J. L. Kelly, ``A New Interpretation of Information Rate," {\it Bell System Technical Journal,}~vol.~35.4,~pp. 917--926,~1956.	
		
		\bibitem{Hakansson_1971}
		N. H. Hakansson, ``On Optimal Myopic Portfolio Policies With and Without Serial Correlation of Yields," {\it Journal of Business,}~vol.~44, \mbox{pp.~324-334},~1971.
		
		%	\bibitem{Cover_Ordentlich_1996}
		%	T. M. Cover and E. Ordentlich, ``Universal Portfolio with Side Information,"~{\it IEEE Transactions on Information Theory,}~IT-42. \mbox{pp.~348-363},~1996.
		
		\bibitem{Merton_1992}
		R. C. Merton, {\it Continuous Time Finance}, Wiley-Blackwell, 1992.
		
		\bibitem{Thorp_2006}
		E. O. Thorp, ``The Kelly Criterion in Blackjack Sports Betting and The Stock Market," {\it Handbook of Asset and Liability Management: Theory and Methodology,}~vol. 1,~pp. 385--428, Elsevier Science,~2006.
		
		\bibitem{Maclean_Thorp_Ziemba_2010}
		L. C. MacLean, E. O. Thorp, and W. T. Ziemba ``Long-term Capital Growth: The Good and Bad Properties of The Kelly and Fractional Kelly Capital Growth Criteria," {\it Quantitative Finance,}~vol.~10,~\mbox{pp. 681--687},~2010.
		
		\bibitem{Cover_Thomas_2006}
		T. M. Cover and J. A. Thomas, {\it Elements of Information Theory}, John Wiley and Sons,~2006.
		
		\bibitem{Luenberger_2013}	
		D. G. Luenberger, {\it Investment  Science},~Oxford University Press,~2013.
		
		\bibitem{Hsieh_Barmish_Gubner_2016_CDC}
		C. H. Hsieh, B. R. Barmish, and J. A. Gubner, ``Kelly Betting Can Be Too Conservative," {\em Proceedings of the IEEE Conference on Decision and Control},~pp.~3695-3701,~Las Vegas,~2016.
		
		\bibitem{MacLean_Thorp_Ziemba_2011}
		L. C. MacLean, E. O. Thorp, and W. T. Ziemba, {\it The Kelly Capital Growth Investment Criterion: Theory and Practice,} World Scientific Publishing Company,~2011.
		
		
		\bibitem{Kuhn_Luenberger_2010}
		D. Kuhn and D. G. Luenberger, ``Analysis of the Rebalancing Frequency in Log-Optimal Portfolio Selection," {\it Quantitative Finance,}~vol.~10,~pp. 221--234,~2010.
		
		\bibitem{Das_Kaznachey_Goyal_2014}
		S. R. Das, D. Kaznachey and M. Goyal, ``Computing Optimal Rebalance Frequency for Log-Optimal Portfolios," {\it Quantitative Finance,}~vol.~14,~pp.1489--1502,~2014.
		
		
		
		%	\bibitem{Zhang_2001}
		%	Q. Zhang, ``Stock Trading: An Optimal Selling Rule," {\em SIAM Journal of Control and Optimization},~vol.-40,~\mbox{pp.~64-87},~2001.
		
		%	\bibitem{Barmish_Primbs_2015_StockTradingViaFeedbackControl}	
		%	B. R. Barmish and J. A. Primbs, ``Stock Trading Via Feedback Control,” {\it Encyclopedia of Systems and Control}, T. Samad and J. Baillieul eds.,~2015.	
		
		
		
		%\bibitem{Hsieh_Barmish_2017_IFAC}
		%C. H. Hsieh and B. R. Barmish, ``On Drawdown-Modulated
		%Feedback in Stock Trading,” {\it IFAC-PapersOnline},~vol.~50, \mbox{pp.~952-958}, 2017.
		
		%
		%
		%
		
		%\bibitem{Hsieh_Barmish_2017_CDC}
		%C. H. Hsieh and B. R. Barmish, ``On Inefficiency of Markowitz-Style Investment Strategies When Drawdown is Important," {\em Proceedings of the IEEE Conference on Decision and Control},~pp. 3075-3080,~Melbourne, Australia,~2017.
		
		\bibitem{Hsieh_Barmish_Gubner_2018_ACC}
		C. H. Hsieh, B. R. Barmish, and J. A. Gubner, ``At What Frequency Should the Kelly Bettor Bet," {\em Proceedings of the  American Control Conference},~pp. 5485--5490,~Milwaukee,~2018.
		
		\bibitem{Hsieh_Gubner_Barmish_2018_CDC}
		C. H. Hsieh, J. A. Gubner, and B. R. Barmish,  ``Rebalancing Frequency Considerations for Kelly-Optimal Stock Portfolios in a Control-Theoretic Framework," {\em Proceedings of the  IEEE Conference on Decision and Control}, pp.~5820--5825,~Miami,~2018.
		
		
		\bibitem{Rujeerapaiboon_Barmish_Kuhn_2018}
		N. Rujeerapaiboon, B. R. Barmish and D. Kuhn, ``On Risk Reduction in Kelly Betting Using the Conservative Expected Value," {\em Proceedings of the IEEE Conference on Decision and Control}, pp.~5801--5806,~Miami,~2018.
		
		
		%	\bibitem{Hsieh_Gubner_Barmish_2018_TAC_Survival}  C. H. Hsieh, B. R. Barmish, and J. A. Gubner, ``A Conjecture Involving Positive Solutions of a Simple Scalar Linear Time-Varying State Equation with Delay," \textit{submitted to the IEEE Transactions on Automatic Control}, 2019.
		
		
		%		\bibitem{Hsieh_Barmish_2015_Allerton}
		%	C. H. Hsieh and B. R. Barmish, ``On Kelly Betting: Some Limitations," {\it Proceedings of the Allerton Conference on Communication, Control, and Computing,}~pp. 165-172,~Monticello,~2015.
		
		
		
		
		%	\bibitem{Finkelstein_Whiteley_1981}
		%	M. Finkelstein and R. Whitley, ``Optimal Strategies for Repeated Games," {\it Advanced Applied Probability},~vol. 13,~pp. 415-428,~1981.
		
		
		
		%	\bibitem{Maclean_Ziemba_Blazenko_1992}
		%	L. C. Maclean, W. T. Ziemba and G. Blazenko ``Growth Versus Security in Dynamic Investment Analysis," {\it Management Science,} vol.~38, pp. 1562-1585,~1992.
		%		
		%	\bibitem{Maclean_Ziemba_1999}
		%	L. C. Maclean and W. T. Ziemba ``Growth Versus Security Tradeoffs in Dynamic Investment Analysis," {\it Annals of Operations Research,} vol.~85, pp. 193-227,~1999.
		
		
		
		\bibitem{Davis_Lleo_2010}
		M. Davis and S. Lleo, ``Fractional Kelly Strategies for Benchmarked Asset Management," in L. C. MacLean, E. O. Thorp, and W. T. Ziemba, {\it The Kelly Capital Growth Investment Criterion: Theory and Practice,} World Scientific,~\mbox{pp. 385-407},~2010.
		
		
		%
		%	\bibitem{Ziemba_2015}
		%	W. T. Ziemba, ``Response to Paul A Samuelson Letters and Papers on the Kelly Capital
		%	Growth Investment Strategy,'' {\it Journal of Portfolio Management}, vol.~42,~pp.~153-167,~2015	
		%	
		
		%	
		%	\bibitem{Billingsley_1995}
		%	P. Billingsley, {\it Probability and Measure,} Wiley-Interscience,~1995.
		
		
		%		\bibitem{Gubner_2006}
		%		J. A. Gubner, {\it Probability and Random Processes for Electrical and Computer Engineers,} Cambridge University Press,~2006.
		%	
		
		
		
		%	\bibitem{Boyd_Vandenberghe_2004}
		%	S. Boyd and L. Vandenberghe, {\it Convex Optimization,} Cambridge University Press, 2004.	
		
		%	
		\bibitem{Barmish_Primbs_TAC_2016}
		B. R. Barmish and J. A. Primbs, ``On a New Paradigm for Stock Trading Via a Model-Free Feedback Controller," {\em IEEE Transactions on Automatic Control},~AC-61,~\mbox{pp.~662-676},~2016.
		
	\end{thebibliography}
\end{document}